%% file: main.tex
\documentclass[10pt]{article}
\usepackage[english]{babel}
\usepackage[letterpaper,top=2cm,bottom=2cm,left=2cm,right=2cm,marginparwidth=1.75cm]{geometry}

\usepackage{float}
\usepackage{calrsfs}
\usepackage{amssymb}
\usepackage{amsmath}
\usepackage{amsthm}
\usepackage{subcaption}
\usepackage{multirow}
\usepackage{graphicx}
\usepackage{amsfonts}
\usepackage{mathtools}
\usepackage{mathrsfs}
\usepackage{placeins}
\usepackage{comment}
\usepackage{bbm}
\usepackage{dsfont}
\usepackage[export]{adjustbox}
\usepackage{overpic}

\usepackage{algorithm}
\usepackage{algpseudocode}

\usepackage{booktabs}

\usepackage[mathscr]{euscript}
\usepackage[dvipsnames,x11names]{xcolor}
\usepackage[colorlinks=true, allcolors=OliveGreen]{hyperref}
\usepackage{amsmath}   
\usepackage{tikz}
\usepackage{pgfplots}
\usetikzlibrary{external}
\usepackage{authblk}
\usepackage{array}
\newcolumntype{C}{>{\centering\arraybackslash}p{2.5cm}}
\usepackage{nicematrix}
\usepackage{caption}
\usepackage{subcaption}
\usepackage{cleveref}
\usepackage{tikz}

\usetikzlibrary{patterns,tikzmark}
\usetikzlibrary{matrix,decorations.pathreplacing,calc}
\usepackage{hf-tikz}

\usepackage{csvsimple}
\usepackage{siunitx}
\usepackage{tikz} 
\usepackage{pgfplots}
\pgfplotsset{compat=newest} 
\usepgfplotslibrary{units} 
\usepackage{tikz}
\usetikzlibrary{matrix,decorations.pathreplacing}
\pgfkeys{/pgfplots/legend entry/.code=\addlegendentry{#1}}
\usepackage{matlab-prettifier}
\usepackage{pdfpages}
\usepackage{enumitem}
\usepackage{orcidlink}

\usepackage[export]{adjustbox}
\usepackage{lipsum}                     
\usepackage{xargs}  

\usepackage{xcolor}
\usepackage{listings}

\definecolor{codebg}{RGB}{250,250,250}
\definecolor{commentgreen}{RGB}{92,128,128}
\definecolor{keywordgreen}{RGB}{120,160,80}
\definecolor{stringred}{RGB}{190,40,40}
\definecolor{linenumber}{RGB}{110,110,110}
\definecolor{framegray}{RGB}{210,210,210}

\lstdefinestyle{mystyle}{
    language=Matlab,
    backgroundcolor=\color{codebg},
    basicstyle=\ttfamily\small,
    columns=fullflexible,
    commentstyle=\color{commentgreen}\itshape,
    keywordstyle=\color{keywordgreen},
    stringstyle=\color{stringred},
    numbers=left,
    numberstyle=\scriptsize\color{linenumber},
    stepnumber=1,
    numbersep=18pt,
    xleftmargin=30pt,
    framexleftmargin=28pt,
    showstringspaces=false,
    tabsize=4,
    breaklines=true,
    frame=single,
    framerule=0.5pt,
    rulecolor=\color{framegray},
    keepspaces=true,
    aboveskip=1em,
    belowskip=1em
}

\usepackage[colorinlistoftodos,prependcaption]{todonotes}
\newcommandx{\unsure}[2][1=]{\todo[linecolor=red,backgroundcolor=red!25,bordercolor=red,#1]{#2}}
\newcommandx{\change}[2][1=]{\todo[linecolor=blue,backgroundcolor=blue!25,bordercolor=blue,#1]{#2}}
\newcommandx{\info}[2][1=]{\todo[linecolor=OliveGreen,backgroundcolor=OliveGreen!25,bordercolor=OliveGreen,#1]{#2}}
\newcommandx{\improvement}[2][1=]{\todo[linecolor=Plum,backgroundcolor=Plum!25,bordercolor=Plum,#1]{#2}}
\newcommandx{\thiswillnotshow}[2][1=]{\todo[disable,#1]{#2}}

\DeclareMathAlphabet{\mathcalligra}{T1}{calligra}{m}{n}

\graphicspath{{Images/}}

\tikzset{  font={\fontsize{15pt}{12}\selectfont}}

\title{The lymph 2.0 library: $p$-adaptive algorithms and parallel assembly strategies for polytopal DG methods\footnote{\textbf{Funding}: This work is partially funded by the European Union (ERC SyG, NEMESIS, project number 101115663). Views and opinions expressed are, however, those of the authors only and do not necessarily reflect those of the European Union or the European Research Council Executive Agency. Neither the European Union nor the granting authority can be held responsible for them. CBLS has been funded by the National Recovery and Resilience Plan (NRRP), Mission 4, Component 1 – Investment 3.4 and Investment 4.1, funded by the European Union. 
MC and CBLS acknowledge “INdAM - GNCS Project”, codice CUP E53C25002010001. The present research is part of the activities of the Dipartimento di Eccellenza 2023-2027 grant, funded by MUR. PFA, SP and CBLS are members of INdAM-GNCS.}}

\author[1]{Paola F. Antonietti\orcidlink{0000-0002-2138-3878}}

\author[1]{Mattia Corti\orcidlink{0000-0002-7014-972X}}

\author[1]{Caterina B. Leimer Saglio\orcidlink{0009-0007-7887-919X}}

\author[1]{Stefano Pagani\orcidlink{0000-0002-6662-3433}}

\affil[1]{MOX-Dipartimento di Matematica, Politecnico di Milano, Piazza Leonardo da Vinci 32, Milan, 20133, Italy}

\begin{document}
\maketitle

\begin{abstract}
This work presents a new release of the \texttt{lymph 2.0} library \cite{antonietti_lymph_2025}, an open-source MATLAB framework for high-order discontinuous Galerkin discretizations on general polytopal meshes. The \texttt{lymph 2.0} version is extended to support discretizations with element-wise polynomial approximation degrees, which allows the design of $p$-adaptive strategies based on a posteriori error indicators. In addition, the library introduces a unified assembly framework that abstracts the construction of discrete operators from the underlying physical model, improving code modularity, parallelism, maintainability, and extensibility. Moreover, the proposed approach enables shared-memory parallelism through dedicated parallel tools. Several numerical examples demonstrate the effectiveness of the proposed developments in reducing the computational cost while preserving approximation accuracy.
\end{abstract}

\input{introduction}
\input{commonAssembly}
\input{Examples}

\section{Conclusions}
\label{sec:conclusions}
The new version of the \texttt{lymph} library introduces a set of advances that turn the original framework into a more powerful environment for high-order discontinuous Galerkin simulations on polygonal meshes. This release aims to improve the generality, modularity, parallelism, and efficiency of the solver by combining a novel assembly engine with an adaptive workflow. The combination of a unified assembly driver with enhanced parallel capabilities improves efficiency for large-scale and long-time simulations, while the reorganization of the postprocessing and input/output layers further streamlines the overall workflow. Moreover, the support for element-wise polynomial degrees and the $p$-adaptive strategies introduces flexible tools to tackle multiscale and strongly heterogeneous problems with a tightly controlled computational cost. Altogether, these developments significantly extend the scope of \texttt{lymph} and strengthen its role as a reliable platform for the development and prototyping of high-order discontinuous Galerkin discretizations on polytopal grids.

\section*{Acknowledgments}
We gratefully acknowledge Stefano Bonetti, Michele Botti, Ivan Fumagalli and Ilario Mazzieri for their careful review and testing of the code released alongside this paper. Their feedback was invaluable in improving the quality, reliability, and usability of the software.

\bibliographystyle{abbrv}
\bibliography{sample.bib}

\end{document}

%% file: introduction.tex
\section{Introduction}
This paper discusses recent updates to the \texttt{lymph} \cite{antonietti_lymph_2025} library (\url{https://bitbucket.org/lymph/lymph}), an open-source MATLAB framework for the numerical approximation of coupled multi-physics problems exploiting polytopal discontinuous Galerkin (PolyDG) discretizations \cite{antonietti_hp_2013,cangiani_hp-version_2014,cangiani_hp-version_2017}. 
The code is specifically designed to address both single-physics and multi-physics problems, offering a modular and extensible structure for mesh handling, finite element space construction, and system assembly. The library supports a variety of applications, including elliptic, parabolic, and hyperbolic problems, as well as coupled systems arising in multi-physics contexts \cite{antonietti_polytopal_2025,antonietti_structure_2026,bonetti_unified_2025,botti_polytopal_2025,leimer_saglio_p-adaptive_2025}.
\par
The present release addresses two important improvements of the library. 
First, we extend the library to support $p$-adaptive discretizations \cite{beirao_da_veiga_manzini_mascotto_2019,houston_wihler_hp_2018,leimer_high-order_2026,melenk_wohlmuth_2001}, allowing for element-wise variation of the polynomial degree, for elliptic and (semilinear) parabolic problems, for which we use the a-posteriori error indicator \cite{cangiani_posteriori_2023,houston_posteriori_2002,leimer_saglio_p-adaptive_2025}. The $p$-adaptive framework reduces the total number of degrees of freedom of the system while retaining a high level of accuracy in the solution approximation. 
Second, we develop a unified assembly framework that abstracts the construction of discrete operators from the specific physical model, enabling a more modular, reusable, and extensible implementation of multi-physics solvers. These improvements significantly simplify the integration of new physical models within the library. In addition, the new assembly algorithm is explicitly designed to exploit shared-memory parallelism, allowing local contributions to be computed concurrently and thereby improving scalability on modern multicore architectures.
Finally, we release two new physics, namely, the Fisher-Kolmogorov equation \cite{fisher_wave_1937,kolmogorov_study_1937} and the monodomain equation coupled with the FitzHugh-Nagumo ionic model \cite{fitzhugh_mathematical_1961,nagumo_active_1962}, that represent two possible examples of nonlinear partial differential equations treatment within the framework of the novel release of the \texttt{lymph 2.0} library.
\par
\begin{figure}[h!]
 \centering
    \includegraphics[width=\textwidth]{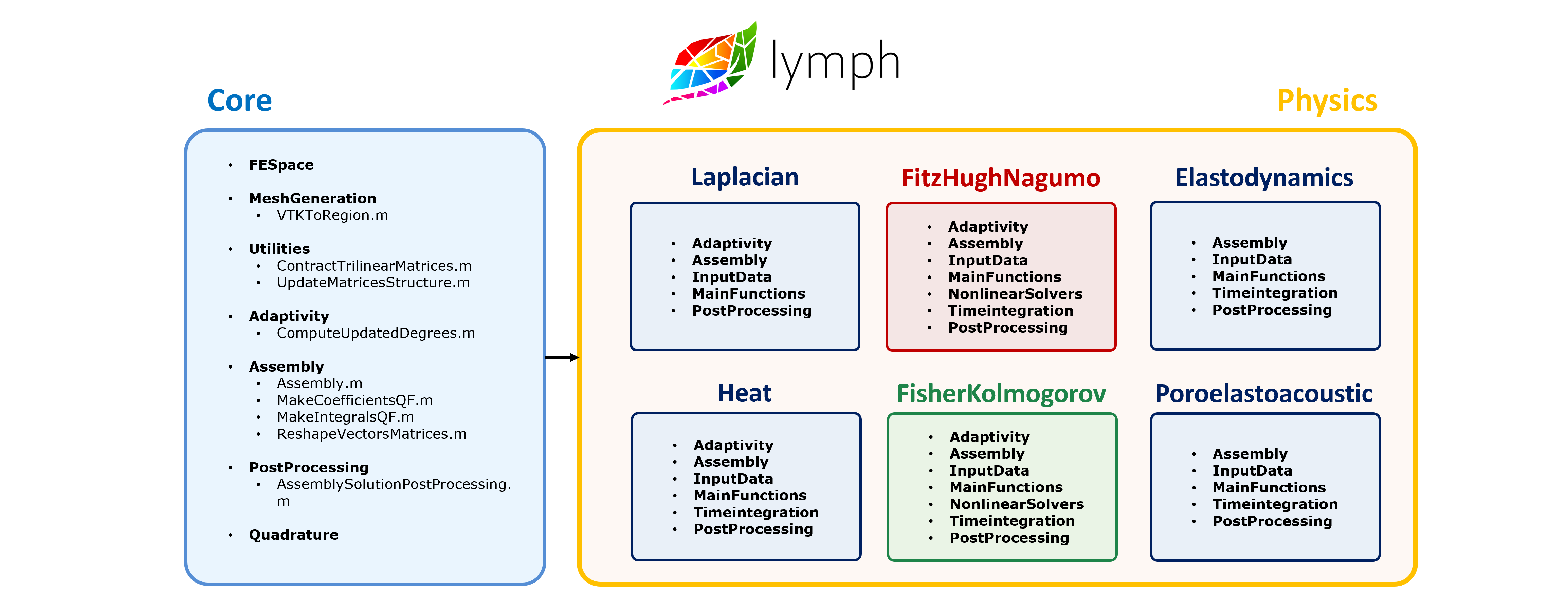}
      \caption{\texttt{lymph 2.0} code structure: highlights of the novel modules and logo (top center).}
  \label{fig:lymph}
  \end{figure}

These new developments significantly improve the flexibility and scalability of the library. The structure of \texttt{lymph 2.0} library with the highlighted modifications is reported in Figure \ref{fig:lymph}.
The paper is organized as follows. Section~\ref{sec:padaptivity} presents the implementation of discretizations with element-wise polynomial approximation degrees, including, for selected physical models, $p$-adaptive strategies. Section~\ref{sec:impro} we resume the main improvements of the library, while Section~\ref{sec:assembly} introduces the new unified assembly framework. Section~\ref{sec:examples} to demonstrate the effectiveness of the proposed features. Finally, conclusions are drawn in Section~\ref{sec:conclusions}.

%% file: commonAssembly.tex

\section{The $p$-adaptive framework}
\label{sec:padaptivity}
One of the main advantages of PolyDG methods is the possibility of using element-wise polynomial approximation orders, thus naturally supporting $p$-adaptive discretization strategies. Indeed, let $p_K\ge 1$ denote the polynomial degree associated with the element $K \in \mathcal{T}_h$. The discrete approximation space is then defined as
\begin{equation*}
V_h^{p} = \left\{ v_h \in L^2(\Omega) : v_h|_K \in \mathbb{P}^{p_K}(K) \quad \forall K \in \mathcal{T}_h \right\}.
\end{equation*}
This setting is especially attractive in multiscale and heterogeneous problems, where different regions of the domain may require different levels of approximation accuracy. In such situations, locally increasing the polynomial degree only where needed may lead to a more efficient use of the degrees of freedom than a globally uniform discretization.

\par
The original \texttt{lymph} framework adopted a discrete finite element space formulated in terms of a global polynomial degree stored in the variable \texttt{femregion.degree}. In this release, we first implement an element-wise discretization degree and a corresponding local number of basis functions, \texttt{femregion.nbases}, for each element. Moreover, we introduce a complete workflow for $p$-adaptivity for the Laplace equation and for semilinear parabolic problems. The routines include a posteriori error indicator computation and efficient updates of the corresponding matrices and finite element structures. From a theoretical point of view, this work is based on the a posteriori error indicator for the polyDG discretization of the Laplacian operator in \cite{cangiani_posteriori_2023} and on the algorithm proposed by \cite{leimer_saglio_p-adaptive_2025} for $p$-adaptivity for semilinear parabolic problems. The latter equations are discretized in time in $(0,T)$ by means of finite difference time stepping $t_0 = t^{(0)} \le t^{(1)} \le ... \le t^{(N)} = T$. A key aspect of the implementation is that the adaptive workflow is entirely driven by the local evaluation of a suitable (residual-based) indicator $\tau_K^{(n)}$, for any $K \in \mathcal{T}_h$ at time $t^{(n)}$. The routine \texttt{ComputeUpdatedDegrees} evaluates the local indicator values and computes a new polynomial distribution according to the adaptive law
\begin{equation*}
   p_K^{\,\mathrm{new}} = f_\mathrm{adapt}\!\left(\frac{\tau_K^{(n)}}{\tau_{\mathrm{threshold}}}\right), \quad K \in \mathcal{T}_h, \; n\in\{0,\dots,N\},
\end{equation*}
where $f_\mathrm{adapt}$ is a user-defined adaptivity function set in \texttt{Data.AdaptFunc}, and $\tau_{\mathrm{threshold}}$ is automatically determined during the first adaptive iteration through a $k$-means clustering procedure applied to the indicator distribution. For a detailed analysis of the algorithm, we refer to \cite{leimer_saglio_p-adaptive_2025}.
\par
To save computational time, the adaptive indicator is not computed for all elements $K$ by default. Indeed, if the solution at a given time snapshot is nearly stationary on the element $K \in \mathcal{T}_h$, the computation of $\tau^{(n)}_K$ is skipped on the element $K$. More precisely, let $U^{(n-1)}$ and $U^{(n)}$ be the expansion coefficients vectors of the time snapshots $t^{(n-1)}$ and $t^{(n)}$, respectively, namely $u_h^{(n)}(\boldsymbol{x}) = \sum_{j=1}^{N_K} U^{(n)}_j \varphi_j(\boldsymbol{x})$, where $\varphi_j(\boldsymbol{x})$ is the set of local basis functions. This choice is driven by the variation of the Euclidean norm $\|U^{(n)} - U^{(n-1)}\|$ at time $t^{(n)}$. The indicator on element $K$ is computed whenever this norm exceeds a (user-defined) tolerance, set by default to $10^{-6}$ in the function \texttt{ComputeAdaptiveIndicatorPhysicsName} for each physics. For the remaining elements, the previously computed indicator values are used, namely $\tau_K^{(n)} = \tau_K^{(n-1)}$.
\par
As a consequence, the $p$-adaptive routine dynamically redistributes the approximation order over the mesh according to the information on the local approximation quality. The implementation also exploits the hierarchical structure of the polynomial basis \cite{antonietti_lymph_2025}. Indeed, if the elementwise polynomial degree is reduced from one adaptive step to the next, the code automatically extract the new matrices associated with volume integrals from the previous higher-order representation, thereby accelerating the reassembly time at each adaptive step.
\paragraph{Example of $p$-adaptive routine: the heat equation.}
Starting from an initial discretization associated with a \texttt{femregion} structure, the \texttt{AdaptivityCycle} routine for each physics performs a sequence of \texttt{Data.AdaptIts} adaptive iterations.
\begin{lstlisting}[style=Matlab-editor, caption={Adaptive cycle for the heat problem.},label={lst:adaptivity-cycle-heat},basicstyle=\footnotesize\ttfamily]
%% Physics/Heat/Adaptivity/AdaptivityCycle.m
    
    AdaptCount = 0;
    while AdaptCount < Data.AdaptIts
        ... % Indicator computation
        [Solution.Indicator] = ComputeAdaptiveIndicatorHeat(.);
        ... % Update the polynomial degrees
        [Data, femregion] = ComputeUpdatedDegrees(.);
        ... % Project the solution onto the new femregion
        [Solution] = ProjectSolutionAssemblyHeat(.);
        ... % Reassemble the matrices
        [Matrices] = MatrixAssemblyHeat(.);
        ...
    end
    
\end{lstlisting}
As shown in Listing~\ref{lst:adaptivity-cycle-heat}, the a posteriori indicator is first computed through the code \texttt{ComputeAdaptiveIndicatorHeat} at each time iteration. The implementation of the indicator is organized within an \texttt{Adaptivity} module in each physics and relies on a sequence of dedicated routines based on the novel common \texttt{Assembly} function. In particular, the adaptive indicator computation involves the following steps:
    \begin{itemize}[itemsep=0.14em]
        \item \texttt{IndicatorPreallocation} initializes the data structures required for the local indicators.
        \item \texttt{ResidualIndicatorAssembly} assembles the cell residual contributions.
        \item \texttt{FacesIndicatorsAssembly} computes the interface and jump terms.
        \item \texttt{FinalIndicator} combines all local contributions into the element-wise adaptive indicator $\tau_K^{(n)}$.
    \end{itemize}
Based on these indicators, the routine \texttt{ComputeUpdatedDegrees} updates the local polynomial distribution over the mesh, producing a new \texttt{femregion} structure with element-wise polynomial degrees adapted to the current solution features. Once the approximation spaces have been modified, the previously computed solution is projected onto the new discrete space by means of the function \texttt{ProjectSolutionAssemblyHeat}. This projection provides a consistent initialization for the subsequent adaptive iteration. The matrices associated with the updated discretization are then reassembled through \texttt{MatrixAssemblyHeat}, accounting for the new local polynomial configuration. The procedure is repeated until the prescribed number of adaptive iterations is reached, progressively refining the distribution of local polynomial degrees and improving the overall approximation efficiency.

\section{An unified and parallel assembly framework}
\label{sec:assembly}
The first release of \texttt{lymph} \cite{antonietti_lymph_2025} already provided efficient strategies for the evaluation of element and face contributions in matrix and vector assembly routines. In particular, the library supported both quadrature-free (QF) \cite{antonietti_fast_2018} and subtriangulation (ST) approaches for polygonal meshes. In this release, we introduce a unified assembly routine that acts as a common driver for the construction of problem-dependent matrices and vectors. The new implementation decouples the general logic of the assembly process from the specific integrals associated with each physical model, resulting in a modular and reusable framework. Moreover, the implementation leverages MATLAB’s \texttt{parfor} loops to assemble local matrices in parallel, improving performance in large-scale computations.
\par
The introduction of a common assembly driver provides several advantages. First, it reduces code duplication by concentrating all quadrature handling, local indexing, and sparse reconstruction in a single routine. As a consequence, it improves maintainability and the improvements to all physical models immediately. Finally, it enhances extensibility, incorporating new physics by defining only the local assembly kernels, without altering the global workflow. Finally, this assembly function is at the basis of both the matrices and forcing terms assembly, as well as the adaptivity indicators and the error computations. In particular, the latter in the first version of the library were computed using the projected solutions instead of the exact ones, potentially leading to approximation issues whenever the underlying solution exhibits sharp fronts.
\begin{lstlisting}[style=Matlab-editor,caption={Common assembly.},label={lst:commonassembly},basicstyle=\footnotesize\ttfamily]
%% Core/Assembly/Assembly.m
    
function [Matrices] = Assembly(Data, neighbor, femregion, AssembInfo, Funcs)
    ... % Preallocation of local matrices for a single element and repetition for all the mesh ones
    Matrices = Funcs.Preallocation(GenMatrices);
    Matrices = repmat({Matrices}, femregion.nel, 1);
    ...
    parfor ie = 1:femregion.nel 
        ... 
        if AssembInfo.assemblyvolumes
            if AssembInfo.ass_vol_vec(ie)
                ... % Assembly the volume
                switch AssembInfo.quadrature
                    case 'QF'
                        ...
                        [Matrices{ie}.Volume]   = Funcs.VolumeAssemblyQF(.);
                        [Matrices{ie}.Volume3L] = Funcs.Volume3LAssemblyQF(.);
                    case 'ST'
                        ...
                        [Matrices{ie}.Volume]   = Funcs.VolumeAssemblyST(.);
                        [Matrices{ie}.Volume3L] = Funcs.Volume3LAssemblyST(.);
                end 
            else
                ... % Reuse and eventual cut of previously assembled matrices
            end
        end
        
        if AssembInfo.assemblyfaces
            ... % Loop over faces
            for iedg = 1 : neighbor.nedges(ie)
                ... % Computation of the matrices related to face integrals
                [Matrices{ie}.Faces] = Funcs.FacesAssembly(.);
            end
        end
    end
    ... % Construction of final global matrices
    [Matrices] = Funcs.FinalMatrices(Matrices);
end
\end{lstlisting}
\paragraph{General structure of the \texttt{Assembly} routine.}
The new framework is built around a single \texttt{Assembly} function (whose structure is reported in Listing~\ref{lst:commonassembly} that has as inputs, besides the mesh and finite element data structures, two problem-dependent ones: \texttt{AssembInfo} and \texttt{Funcs}. The first collects all the information needed to control the assembly workflow. This includes the choice of the quadrature strategy for volume terms, logical flags determining which types of integrals must be assembled (i.e., volume, external/internal face, trilinear terms), and additional options related to basis-function derivatives and adaptive updates. These flags allow the code to skip unnecessary computations in specific assembly calls, thereby preserving the original efficiency. Additionally, the structure contains the current time or solution values for the assembly of nonlinear or forcing terms.
\par
The \texttt{Funcs} structure instead contains the problem-specific routines that define the contributions for the specific physical problem. Namely, the structure requires a handle for matrix preallocation, volume, face, and trilinear form assembly, as well as final construction of the global matrices. As a consequence, the same assembly engine can be reused for different applications by simply changing the local assembly kernels collected in \texttt{Funcs}.
\par
The unified routine preserves at the volume level the support of both QF and ST strategies, consistently with the original release \cite{antonietti_lymph_2025}. The two different implementations for the specific physics need to be separately provided through the handles associated with \texttt{Funcs.VolumeAssemblyQF} and \texttt{Funcs.VolumeAssemblyST}, respectively. In any case, only the strategy chosen by the flag \texttt{AssembInfo.quadrature} will be used in practice. As an additional advantage, in the subtriangulation assembly of the volume terms, we eliminated an unnecessary loop over the resulting triangles, improving the overall computational time. On the contrary, the face terms are always computed using a subtriangulation of the segment, without making use of QF strategies. For this reason, only a function \texttt{Funcs.FacesAssembly} needs to be constructed. As in the original release, the assembly strategy is based on the construction of local matrices and index structures used to map local degrees of freedom to the global algebraic system. This process can now be performed in an element-wise parallel approach by using the \texttt{parfor} construct to distribute the assembly strategy across multiple cores, by setting \texttt{Setup.isParallel} equal to 1. Once all local contributions have been computed, the data are reshaped, filtered to remove unused entries, and passed to the final reconstruction stage. In the common assembly framework, only the assembly of local matrices is under the control of the single physics, while the construction of the global ones is automated and contained in the \texttt{ReshapeVectorMatrices} function. The final matrices associated with the discretization method are then constructed through the problem-specific \texttt{Funcs.FinalMatrices} routine.
\paragraph{Integration with $p$-local workflows.} 
The \texttt{Assembly} function has been adapted to account for element-wise varying polynomial degrees. The design extends to adaptive settings, where local polynomial degrees may be updated at different stages. To avoid unnecessary computations, in the $p$-adaptivity steps, the assembly process considers only the elements contained in the vectors \texttt{AssembInfo.ass\_vol\_vec} and \texttt{AssembInfo.ass\_face\_vec}, based on the local degree changes. On the contrary, for the remaining elements the assembly process can be skipped by reusing the previously computed matrices or by extracting them in case of reduction of the polynomial degree $p$.
\paragraph{Nonlinear PDEs.} Another important feature of the new assembly routine is the possibility of handling, within the same framework, nonlinear contributions as well.
Concerning general forms of nonlinearities, the assembly is typically called inside the nonlinear solver. To handle the presence of previous iteration solutions, we manage the passage of the information within the \texttt{AssembInfo} structure. The assembly of the specific structures in this case is performed within the classical \texttt{VolumeAssembly} and \texttt{FacesAssembly} functions. On the contrary, we provide a separate treatment for trilinear contributions. Indeed, the quadratic reaction terms of semilinear equations are controlled by the logical flag \texttt{assemblytrilinearforms}, which activates \texttt{VolumeAssembly3L} for the construction of three-dimensional tensors $[\mathbf{M}]_{ijk} = \int_\Omega \varphi_i(\boldsymbol{x}) \varphi_j(\boldsymbol{x}) \varphi_k(\boldsymbol{x}) \mathrm{d}\boldsymbol{x}$. 
For these terms (e.g., in the Fisher-Kolmogorov equation), the use of the block-diagonal structure of the volume matrices allows contracting the resulting tensor into a matrix using the tool \texttt{ContractTrilinearMatrix}, which is more efficient than reassembling it at every iteration.

\paragraph{Multiphysics problems.} This unified perspective is particularly advantageous for multiphysics problems. Different operators, potentially associated with distinct fields or submodels, can be assembled by combining appropriate local routines while relying on a common high-level driver. As a result, the implementation becomes more robust, easier to maintain, and naturally extensible to new coupled formulations. Moreover, the approach enables the assembly of all physical contributions within a single loop, avoiding repeated preallocation and quadrature steps.
\section{Additional improvements}
\label{sec:impro}
In addition to multiple bug fixes and performance improvements, the new release includes the following novel improvements:
\begin{itemize}
\item \textbf{New models released in Physics:} In this release, we introduce the Fisher-Kolmogorov equation and monodomain equation coupled with the FitzHugh-Nagumo ionic model as examples of semilinear PDEs stemming, for example, in the fields of biological population dynamics and electrophysiology.
\item \textbf{Improved postprocessing assembly:} As done for the matrices, the postprocessing is now based on a common function \texttt{AssemblySolutionPostProcessing} that minimizes code duplication in loops.
\item \textbf{Distributed matrices in parallel:} In case of large computations, it can be useful to distribute the matrices between different cores to save computational time. An example of the use of MATLAB's \texttt{distribute} function is provided in the \texttt{Elastodynamics} module.
\item \textbf{Background postprocessing in parallel framework:} Whenever the parallel flag \texttt{Setup.isParallel} is activated, the library automatically makes use of the \texttt{parfeval} function for postprocessing solutions, performing it in the background without slowing the evaluation of subsequent time steps.
\item \textbf{Input VTK \cite{vtkBook} meshes:} The library can now read meshes in the \texttt{.vtk} format. The function \texttt{VTKtoRegion} automatically handles the generation of the \texttt{region} and \texttt{neighbor} structures and saves them in a \texttt{.mat} file for future use.
\end{itemize}

%% file: Examples.tex
\section{Numerical examples}
\label{sec:examples}
In the following, we show some paradigmatic examples to illustrate the new capabilities of \texttt{lymph} in handling both steady-state and time-dependent, as well as scalar and vector-valued PDEs.
\subsection{The Poisson problem: $p$-adaptive solver}
We start by considering the Poisson problem in a polygonal domain $\Omega \subset \mathbb{R}^2$:
\begin{equation}
\label{eq:laplacian}
\begin{aligned}
-\nabla \cdot (\mu \nabla u) &\, = f, && \qquad \mathrm{in}\;\Omega, \\
u &\,= g, && \qquad \mathrm{on}\;\partial \Omega.
\end{aligned}    
\end{equation}
The numerical discretization of problem \eqref{eq:laplacian} is performed by means of the PolyDG method described in \cite{antonietti_lymph_2025}. As a first test case, we consider the domain $\Omega=(0,1)^2$ with homogeneous diffusion coefficient $\mu = 1$. The exact solution is $u(x,y)=\tanh\!\left(-20(x^2+y^2-0.8)\right)$, and the boundary condition and the forcing term are computed accordingly. To solve this problem, we use the routines contained in the \texttt{Laplacian} module, and the data are set in the file \texttt{DataComparisonpAdaptive.m}. 
We focus on the comparison between adaptive and uniform simulations using an adaptive procedure to dynamically update the local polynomial degree distribution up to $p_{\max}=10$. We exploit different element meshes ($N_\mathrm{el} = 180, 300, 800, 1500, 2200, 3000$). The $p$-adaptive indicator used in the simulation is characterized by the sum of a residual component and one associated with the jumps of the numerical solution and its gradient across the mesh edges (see \cite{cangiani_posteriori_2023} for details).
\begin{figure}[h!]
 \centering
  \begin{subfigure}[t]{0.45\textwidth}
  \centering
  \raisebox{5mm}{
\includegraphics[width=1.1\textwidth]{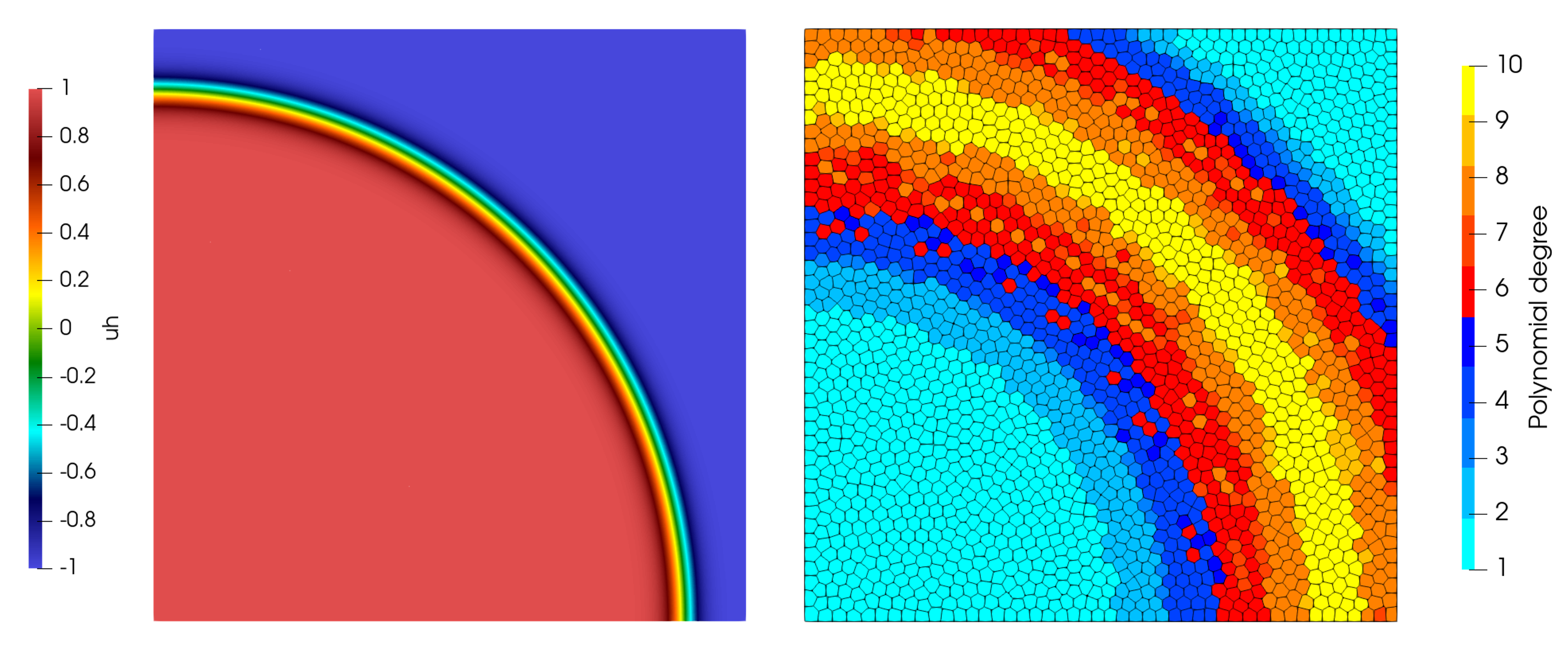}
    }
  \end{subfigure}\hspace{3em}
\begin{subfigure}[t]{0.4\textwidth}
  \centering
    \includegraphics[width=\textwidth]{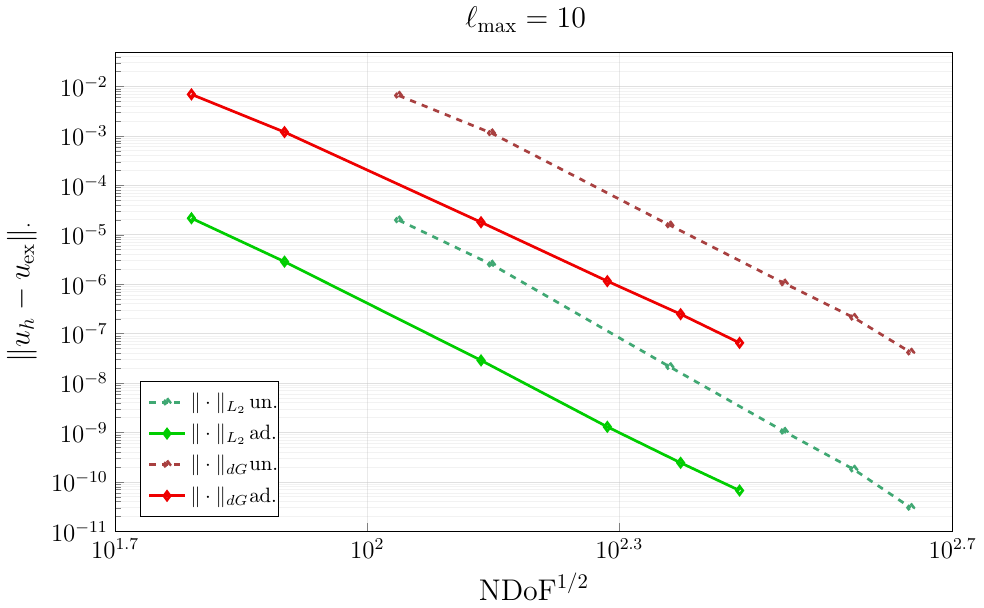}
  \end{subfigure}
      \caption{Poisson problem. Left: Numerical solution $u_h$ and polynomial degree distribution with $p_\mathrm{max}=10$. Right: Corresponding computer errors in the $\|u-u_h\|_{dG}$ and $\|u-u_h\|_{L^2(\Omega)}$ as a function of the total number of degrees of freedom of the system. The adaptive approach (solid line) achieves a reduction of approximately $38\%$--$43\%$ degrees of freedom compared to the uniform refinement (dashed line).}
  \label{fig::laplacian1}
\end{figure}
\par
Figure~\ref{fig::laplacian1} illustrates the behavior of the $p$-adaptive strategy compared with the uniform one. On the left, we show the computed numerical solution together with the computed element-wise polynomial degree distribution obtained with the adaptive algorithm compared to the uniform refinement (dashed line), which is obtained by exploiting a uniform polynomial degree. On the right, we show the computed errors as a function of the total number of degrees of freedom (NDoF) of the system in the $L^2$ \textit{dG}-norms, where the latter is defined as: 
\begin{equation}
\|u\|_{dG}^2 = \|\mu\nabla_h u\|_{\mathbf{L}^2(\Omega)}^2+\|\sqrt{\eta}[\![u]\!]\|^2_{L^2(\mathcal{F})}.
\end{equation}
Figure~\ref{fig::laplacian1} shows that high polynomial degrees are automatically assigned along the transition layer.
The convergence behavior is preserved. Here, the same level of accuracy is achieved with significantly fewer NDoF. This confirms that the $p$-adaptive approach is able to maintain the approximation properties of the uniform high-order discretization while substantially reducing computational cost.

\subsection{Time-dependent problems}
As a second example, we consider a general semilinear parabolic problem. Given an open, bounded, polygonal domain $\Omega \subset \mathbb{R}^d$, $(d=2,3)$ and a final time $T>0$, we consider the solution $u: \Omega \times [0,T] \rightarrow \mathbb{R}$. The model reads as follows: for any time $t \in (0,T]$, find $u=u(\boldsymbol{x},t)$ such that:
\begin{equation}
\label{eq:general}
\frac{\partial u}{\partial t} - \nabla \cdot (\mu \nabla u) + f(u) = g_u, \qquad \mathrm{in}\;\Omega,
\end{equation}
complemented with suitable boundary and initial conditions $u(0) = u^0$. Moreover, $f(u)$ is a possibly nonlinear function of the solution, and $g_u$ is the forcing term.
\subsubsection{The heat equation}
First, we consider the heat equation in \eqref{eq:general} setting $f=0$. Moreover, we consider $\Omega=(-1,1)^2$ and $\mu=1$. The exact solution is given by
\[
u(x,y,t) = \exp\left(\frac{1 - \exp\left(\frac{x - t y}{\varepsilon}\right)}{1 - \exp\left(-\frac{1}{\varepsilon}\right)}\right),
\]
where $\varepsilon=8e-2$ is the parameter controlling the thickness of the resulting layer in the solution. The forcing term and the Dirichlet boundary conditions are set accordingly. To solve this problem, we use the routines contained in the \texttt{Heat} module using the data in \texttt{DataTestCasepAdaptive.m}. For this simulation, we impose $T = 3$ and $\Delta t = 10^{-2}$, with implicit Euler scheme discretization. The mesh is composed of 1000 polygonal elements ($h = 0.121$) with $p_{\max}=5$. The numerical simulations have been performed as a serial job using the \texttt{Kami} cluster (40 computing nodes configured as follows: \textbf{CPU} 2$\times$ AMD EPYC 7413 24-Core Processor, \textbf{RAM} 512 GB) at the Department of Mathematics, Politecnico di Milano.

\begin{table}[ht]
\centering
\resizebox{\textwidth}{!}{%
\begin{tabular}{|c|c c c c c c c|c|}
\hline
  & Matrices
  & Indicator 
  & RHS 
  & Degree Distribution 
  & Numerical Solution 
  & Algebraic System 
  & Forcing Term 
  & \textbf{Total} \\
  & Reassembly
  & Construction 
  & Assembly
  & Computation 
  & Projection 
  & Solving
  & Computation 
  & \textbf{Time} \\
\hline

\textbf{Uniform} 
& - 
& - 
& 0.54 s
& - 
& - 
& 193.59 s
& 143.93 s
& \textbf{338.06 s} \\

{\textit{\% of total}}
& - 
& - 
& \textit{0.16\%}
& - 
& - 
& \textit{57.27\%}
& \textit{42.57\%}
&  \\

\hline
\textbf{Adaptive} 
& 27.54 s
& 30.47 s
& 0.23 s
& 0.29 s
& 11.47 s
& \phantom{0}62.07 s
& \phantom{0}83.93 s
& \textbf{216.00 s} \\

{\textit{\% of total}}
& {12.75\%}
& \textit{14.11\%}
& \textit{0.11\%}
& \textit{0.13\%}
& \textit{5.31\%}
& \textit{28.74\%}
& \textit{38.86\%}
&  \\

\hline
\end{tabular}%
}
\caption{Heat equation. Computational times for uniform and adaptive strategies. Percentages are computed with respect to the total time of each strategy.}
\label{tab:timings}
\end{table}

\par
\begin{figure}[h!]
 \centering
\includegraphics[width=\textwidth]{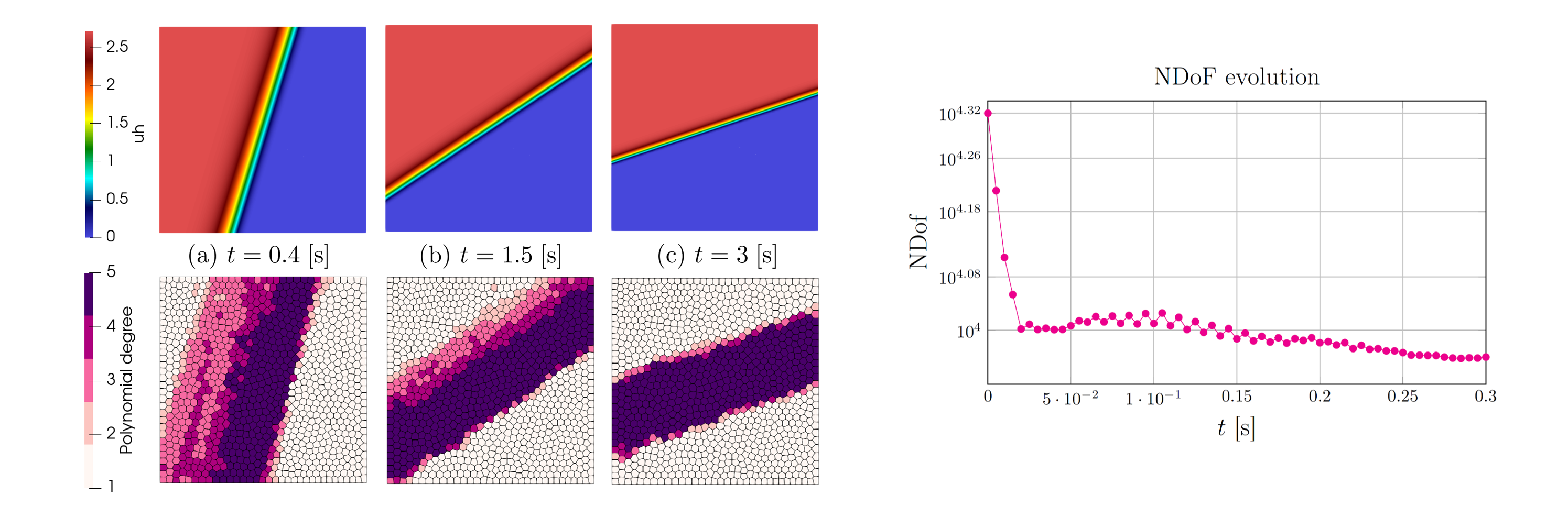}
      \caption{Heat equation. Evolution of the numerical solution and corresponding polynomial degree distribution at different time instants \(t=0.4\,\mathrm{s}\), \(t=1.5\,\mathrm{s}\), and \(t=3\,\mathrm{s}\).}
  \label{fig:heat1}
\end{figure}

Table~\ref{tab:timings} reports the computational times for the uniform and $p$-adaptive strategies for this simulation. Although the adaptive approach introduces additional costs related to matrix reassembly, indicator construction, and degree distribution within the $p$-adaptive routine, it significantly reduces the time required to solve the linear system and to compute the forcing term of the problem. As a result, the overall computational time is reduced from $338.06\,\mathrm{s}$ to $216.00\,\mathrm{s}$, corresponding to a gain of approximately $36.11\%$. Figure~\ref{fig:heat1} shows the evolution of the total NDoFs during the simulation. The evolution highlights the effect of the adaptive strategy on the total number of degrees of freedom of the system. In particular, the adaptive approach provides a reduction of approximately $54\%$ in the number of degrees of freedom with respect to the initial configuration with $p_{\max} = 5$.

\subsubsection{The Fisher-Kologorov model}
In this section, we consider the Fisher-Kolmogorov equation obtained from problem \eqref{eq:general} by choosing $f(u) = - u(1-u)$. This model is widely employed in biological wave propagation and represents an example of a nonlinear parabolic equation exhibiting traveling fronts and sharp transition layers. 
\par
\begin{figure}[h!]
 \centering
\includegraphics[width=\textwidth]{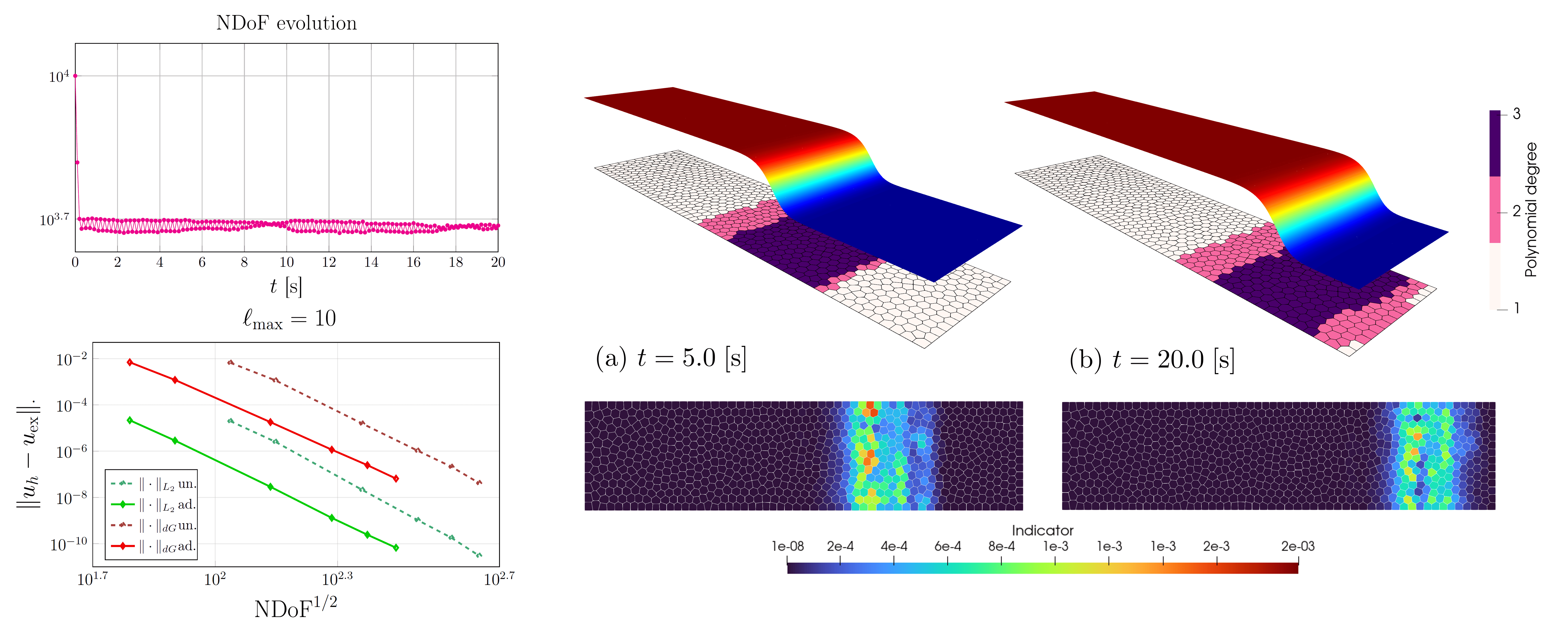}
\caption{Fisher-Kolmogorov model. Left:  temporal evolution of the total number of NDoFs during the
$p$-adaptive simulation and convergence history in the $L_2$ and $dG$-norms with respect to $\mathrm{NDoF}^{1/2}$. Right: numerical solution and spatial distribution of the local polynomial degree and a-posteriori indicator for each cell.}
  \label{fig:FK_solution}
\end{figure}

For details about the discretization of the problem, we refer to \cite{corti_discontinuous_2025}. The computational domain is $\Omega = [-1,3]\times[0,1]$, with final time $T=4$. The coefficients are $\mu = 10^{-3}$, the polynomial degree is $p_{\max}=3$, and the time step is $\Delta t=10^{-2}$. The exact solution of the problem is set to $u(x,y,t)=0.25[1+\tanh(8-\sqrt{(24\mu)^{-1}}\,x)]^2$.
We exploit homogeneous Neumann boundary conditions as in \cite{antonietti_structure_2026}, and the nonlinear term is treated by Picard iterations with tolerance $10^{-10}$ and maximum number of iterations equal to $1000$. Finally, adaptive polynomial refinement is employed. The adaptive algorithm is performed every 10 time steps, and at most two adaptation iterations are allowed. To solve this problem, we use the functions contained in \texttt{FisherKolmogorov}, and we set the data in \texttt{DataWavesFKPP.m}. Figure~\ref{fig:FK_solution} (right) shows the numerical solution at different time instants $t=5.0$ and $t=20\;s$, the corresponding polynomial degree distribution obtained through the adaptive procedure, and the distribution of the complete a posteriori indicator $\tau^{(n)}_K$. As expected, the indicator attains its largest values in correspondence with the sharp transition layer of the traveling wave, while remaining small in the smooth regions of the solution, so that the highest approximation orders are concentrated along the moving front. 
In Figure~\ref{fig:FK_solution} (left), we report a convergence analysis for the uniform refinement with final time $T = 1$ and time step $\Delta t = 0.01$. Different mesh resolutions are considered, with a polynomial degree $p_{\max}=3$. The convergence behavior of the uniformly enriched discretization is compared with that obtained through the adaptive strategy. As expected, the adaptive approach preserves the same level of accuracy achieved by the uniform approximation while significantly reducing the NDoFs.
%
%
\subsubsection{The monodomain coupled with FizHugh-Nagumo model}
In this section, we consider the monodomain equation coupled with FitzHugh-Nagumo ionic model, which provides a simplified description of the electrophysiological behavior of biological tissues \cite{fitzhugh_mathematical_1961, nagumo_active_1962}.
Given an open, bounded, polygonal domain $\Omega \subset \mathbb{R}^d$, $(d=2,3)$ and a final time $T>0$, we consider the solution $(u,{w})$ with $u: \Omega \times [0,T] \rightarrow \mathbb{R}$, and ${w} = {w}(\boldsymbol{x},t)$ with ${w}: \Omega \times [0,T] \rightarrow \mathbb{R}$. The model reads as follows: for any time $t \in (0,T]$, find $u=u(\boldsymbol{x},t)$ and $w=w(\boldsymbol{x},t)$ such that:
\begin{equation}
\label{eq:FHN}
\begin{aligned}
C \chi_m\frac{\partial u}{\partial t} - \nabla \cdot (\mu \nabla u) + \chi_m ku(u-a)(u-1) &\, = g_u, && \qquad \mathrm{in}\;\Omega, \\
\frac{\partial w}{\partial t} + \varepsilon(\beta u - \gamma w) &\,= g_w, && \qquad \mathrm{in}\;\Omega.
\end{aligned}    
\end{equation}
complemented with suitable boundary and initial conditions $u(0) = u^0$ and $w(0) = {w}^0$. With parameters $a \in (0,1)$, $\varepsilon \ll 1$, and $\gamma > 0$, controlling the excitability and recovery dynamics of the system. In this test case, we consider $\Omega=(0,2.5)^2$ with homogeneous Neumann boundary conditions. We set the model parameters as $\mu = 10^{-4}$, $\kappa = \gamma = 1$, $\varepsilon = 0.003$, $\beta = 0.5$, and $a = 0.1$. We adopt a Crank–Nicolson time-stepping scheme with final time $T=3000$ and time step $\Delta t = 0.5$. We consider a Cartesian mesh of $1\,024$ elements and employ the $p$-adaptive DG approximation with a polynomial degree $p_{\max}=4$. Adaptive updates are performed every $10$ time steps, with at most $2$ adaptive iterations. The resulting nonlinear algebraic system is solved through Picard iterations with tolerance $10^{-10}$ and maximum number of iterations equal to $1000$. 
The initial conditions are chosen in order to trigger the formation and propagation of spiral-wave dynamics, as follows 
\[
u_0(x,y)=
\begin{cases}
1, & x \leq 1.25 \ \text{and} \ y \leq 1.25,\\
0, & \text{otherwise},
\end{cases}
\qquad w_0(x,y)=\frac{1}{20}\left(1+\tanh\left(\frac{y-1.25}{0.02}\right)\right).
\]
\begin{figure}[!hbtp]
 \centering
    \includegraphics[width=\textwidth]{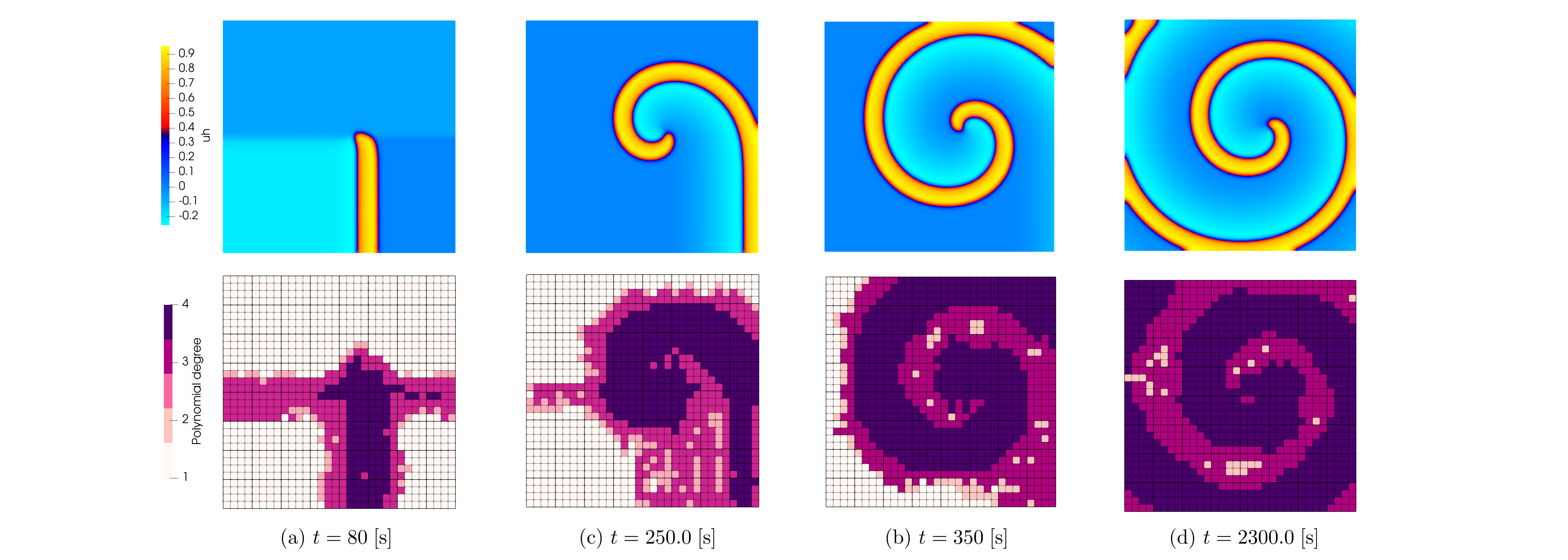}
      \caption{FitzHugh-Nagumo model. Snapshots of the numerical solution \(u_h\) (top row) and corresponding polynomial-degree distribution (bottom row) at different time snapshot $t=80,\;250,\;350,\;2300 \;{s}$.}
  \label{fig::fhn1}
  \end{figure}
\par
To solve this problem, we use the \texttt{FitzHughNagumo} module and the data file \texttt{DataFHNSpiral.m}. Figure~\ref{fig::fhn1} shows the numerical solution of the FitzHugh–Nagumo system at four representative time instants. The first row shows the evolution of the activator variable $u_h$, while the second row reports the corresponding distribution of the local polynomial degree selected by the adaptive strategy. We observe that the highest polynomial degrees are located in proximity to the spiral wave. 
\subsection{The elastodynamic system}
We now consider the linear elastodynamic equation: find the displacement field $\boldsymbol{u} : \Omega \times [0,T] \to \mathbb R^2$ such that
\[
\rho \partial_{tt} \boldsymbol{u} - \nabla \cdot \boldsymbol{\sigma}(\boldsymbol{u}) = \boldsymbol{f}
\qquad \text{in}\;\Omega \times (0,T],
\]
together with proper boundary and initial conditions $\boldsymbol{u}(\cdot,0)=\boldsymbol{u}_0,\; \partial_t \boldsymbol{u}(\cdot,0)=\boldsymbol{u}_0$. The stress tensor is defined by the isotropic linear elastic constitutive relation $\boldsymbol{\sigma}(\boldsymbol{u})=2\mu \boldsymbol{\varepsilon}(\boldsymbol{u})
+\lambda \operatorname{tr}(\boldsymbol{\varepsilon}(\boldsymbol{u})) \mathbf{I}$ , where $\boldsymbol{\varepsilon}(\boldsymbol{u})$ is the symmetric strain tensor, $\mathbf{I}$ denotes the identity tensor, and $\lambda,\mu \in L^\infty(\Omega)$ are the Lamé coefficients. All the details about material parameters, model coefficients, boundary conditions and numerical discretization adopted in this work are reported in \cite[\S5.3]{antonietti_lymph_2025}. The numerical simulations have been performed as a serial job using the \texttt{Nemesis} cluster (\textbf{CPU} 2$\times$ AMD EPYC 9634 84-Core Processor (336 threads), \textbf{RAM} 1.5 TB) at the Department of Mathematics, Politecnico di Milano.

\begin{table}[!hbtp]
\centering
\begin{tabular}{|r|c c c c|}
\hline
& Matrix Assembly (QF) 
& Matrix Assembly (ST) 
& RHS Assembly 
& Linear System Solver  \\
\hline
\textbf{Serial: 1 core\phantom{s}} 
&  9.51 s 
& 10.23 s
&  1.93 s 
& 60.06 s 
\\ (\% time saving)
& \textit{(-2.36\%)} 
& \textit{(-4.57\%)} 
& \textit{(-4.46\%)} 
& \textit{(-1.54\%)} 
\\ \hline
\textbf{Parallel: 4 cores}
& 5.05 s 
& \phantom{0}5.39 s 
& 0.72 s
& 26.20 s 
\\ (\% time saving)
& \textit{(-48.15\%)} 
& \textit{(-49.72\%)} 
& \textit{(-64.36\%)} 
& \textit{(-57.05\%)} 
\\ \hline
\textbf{Parallel: 8 cores}
& 4.31 s 
& \phantom{0}4.58 s 
& 0.47 s
& 22.13 s 
\\ (\% time saving)
&  \textit{(-55.75\%)} 
&  \textit{(-57.28\%)} 
&  \textit{(-76.73\%)} 
&  \textit{(-63.74\%)} 
\\
\hline
\textbf{\texttt{lymph} 1.0 \cite{antonietti_lymph_2025}} 
& 9.74 s 
& 10.72 s 
&  2.02 s 
& 61.00 s  
\\
\hline
\end{tabular}
\caption{Elastodynamics equation. Computational times serial and parallel solvers and comparison with the \texttt{lymph}~1.0 version. With total number of dregrees of freedom of the system equal to 102270.}
\label{tab:times_elasto}
\end{table}

\par
Table~\ref{tab:times_elasto} reports the computational times for the different phases of the simulation, comparing serial and parallel runs with the new assembly framework against the previous \texttt{lymph}~1.0 implementation. The results show that the unified assembly routine benefits significantly from shared-memory parallelism: moving from a serial run to 4 and 8 cores yields a substantial reduction in the matrix assembly times, for both QF and ST strategies, and for the linear system solver time. Finally, the row corresponding to \texttt{lymph}~1.0 highlights that the new release achieves comparable times in serial, confirming the effectiveness of the redesigned driver, with important gains obtained when parallelization is enabled.